\documentclass[11pt,a4paper,reqno]{amsart}

\usepackage{amsmath}
\usepackage{mathptmx}
\usepackage[scaled=.90]{helvet}
\usepackage{courier}

\usepackage{graphicx,parskip}
\newtheorem{theorem}{Theorem}
\newtheorem{lemma}[theorem]{Lemma}
\newtheorem{corollary}[theorem]{Corollary}
\theoremstyle{definition}
\newtheorem{example}[theorem]{Example}

\begin{document}
\title[Constructing interpolating Blaschke products]%
{Constructing interpolating Blaschke products\\with given preimages}
\date{March 31, 2006}
\author{Geir Arne Hjelle}
\urladdr{http://www.math.ntnu.no/\textasciitilde hjelle/}
\email{hjelle@math.ntnu.no}
\address{Geir Arne Hjelle\\Department of Mathematical Sciences\\Norwegian University of Science and Technology\\7491 Trondheim\\Norway}
\thanks{The research is partially supported by a grant from the Research Council of Norway, project \#155060.}

\begin{abstract}
We give a constructive and flexible proof of a result of P.\ Gorkin
and R.\ Mortini concerning a special finite interpolation problem on
the unit circle with interpolating Blaschke products. Our proof also
shows that the result can be generalized to other closed curves than
the unit circle.
\end{abstract}

\maketitle

\section{Introduction}
\label{sec:introduction}
Let $\mathbb D$ denote the unit disk $\{z \colon |z| < 1 \}$.  A
finite Blaschke product of degree $N$ is a function $B \colon
\overline{\mathbb D} \rightarrow \overline{\mathbb D}$ of the form
\begin{equation*}
  B(z) = \lambda \prod_{n=1}^N \frac{z - z_n}{1 - \bar z_n z} ,
\end{equation*}
where $|\lambda| = 1$ and $z_n \in \mathbb D$ for $n = 1, \ldots, N$. A
finite Blaschke product of degree $N$ maps the unit circle $\partial
\mathbb D$ onto itself $N$ times. The uniform separation constant of $B$ is
\begin{equation*}
  \delta(B) = \min_k \prod_{j \neq k}
    \biggl| \frac{z_j - z_k}{1 - \bar z_k z_j} \biggr| .
\end{equation*}
A Blaschke product $B$ is a finite interpolating Blaschke product if
all the zeros $z_n$ of $B$ are simple, that is if $\delta(B) > 0$.  We
will show the following result.

\begin{theorem}
\label{thm:mainthm}
Let $\{\gamma_n \colon n = 1, \ldots, N\}$ be a partition of the unit
circle into a finite number of arcs. For every $C < 1$, there is a
finite Blaschke product $B$ of degree $N$ such that $\delta(B) > C$
and $B(\gamma_n) = \partial \mathbb D$ for $n = 1, \ldots, N$.
\end{theorem}

The result was proved by P.\ Gorkin and R.\ Mortini
in~\cite{Gorkin05}, using a non-constructive argument based on a
theorem of W.\ Jones and S.\ Ruscheweyh~\cite{Jones87} (see
Theorem~\ref{thm:Jones87} below). We will present a constructive and
flexible proof. The proof will also show that the result can be
generalized to other closed curves than the unit circle. We discuss
this extension in Corollary~\ref{cor:extension}.

To prove Theorem~\ref{thm:mainthm} we will describe an iterative
algorithm that constructs a sequence of Blaschke products $(B_k)$
converging to the desired Blaschke product $B$. The algorithm, which
is inspired by the circle packing algorithm of Collins and
Stephenson~\cite{Collins03}, depends on certain monotonicity
relations. The details are found in Section~\ref{sec:algorithm}. The
last section of this paper explores how Theorem~\ref{thm:mainthm}
relates to radial limits and finite interpolation problems.

Theorem~\ref{thm:mainthm} can be formulated as a problem of finite
interpolation. Given $N$ distinct points $\varphi_1, \ldots,
\varphi_N$ on the unit circle $\partial \mathbb D$, there is a finite
Blaschke product, $B$, of degree $N$, such that $\delta(B) > C$ and
$B(\varphi_n) = 1$ for $n = 1, \ldots, N$. If we drop the condition on
$\delta(B)$, the result follows readily from the following result of
Jones and Ruscheweyh~\cite{Jones87}.

\begin{theorem}
\label{thm:Jones87}
Let $\varphi_1, \ldots, \varphi_N \in \partial \mathbb D$ be distinct, and
let $\psi_1, \ldots, \psi_N \in \partial \mathbb D$. Then there exists a
Blaschke product, $B$, of degree at most $N - 1$ satisfying
\begin{equation}
\label{eq:fip}
  B(\varphi_n) = \psi_n, \qquad n = 1, \ldots, N .
\end{equation}
\end{theorem}

In our case we must add one equation to avoid the trivial solution
$B(z) \equiv 1$. We choose $\varphi_{N+1} \in \partial \mathbb D$ such that
$\varphi_{N+1} \neq \varphi_n$ for $n = 1, \ldots, N$, and demand that
\begin{equation*}
  B(\varphi_{N+1}) = \psi_{N+1} \neq 1.
\end{equation*}
Unfortunately, the proof of Jones and Ruscheweyh is non-constructive,
and gives little information about the localization of the zeros of
$B$. Hence, in order to get a constructive proof of
Theorem~\ref{thm:mainthm} a different argument is needed.

\begin{figure}[bt]
\begin{center}
\includegraphics{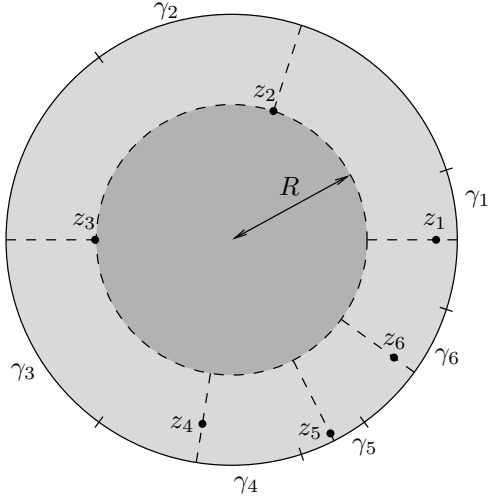}
\caption{The zeros are placed on radii at least a distance $R$ away
from the origin, that is, in the lightly shaded region.}
\label{fig:extraconditions}
\end{center}
\end{figure}

There is a lot of freedom in the problem. In~\cite{Semmler06} G.\
Semmler and E.\ Wegert discuss the uniqueness of solutions of finite
interpolation problems on the form~\eqref{eq:fip} with minimal
degree. Our problem is what they call damaged. In particular, this
means that the minimal degree solution is not unique. The lack of
uniqueness stems from the following informal argument. In order to
place $N$ zeros in $\mathbb D$, we have to decide the value of $2N$ real
variables.  However, the arcs on $\partial \mathbb D$ gives rise to only
$N-1$ equations, as one equation can always be satisfied with a
properly chosen rotation. We use this freedom to add $N+1$ extra
conditions.
\begin{enumerate}
\item[i)] One zero is placed on each of the radii through the
mid-point of each arc.
\item[ii)] All zeros are placed at least a distance $R > 0$ away from
the origin.
\end{enumerate}
See Figure~\ref{fig:extraconditions}. It is not obvious that we can
still solve the problem under these extra conditions. The point is
that i) and ii) guarantee that if a solution $B$ exists, then its
zeros will be simple and $\delta(B)$ will be bigger than some constant
depending on the arcs $\{\gamma_n\}$ and $R$.

Before focusing on the general case, we find the solutions in some
specific examples.

\begin{example}
Assume all arcs have the same length, $\frac{2\pi}{N}$. In this case,
symmetry considerations imply that if the zeros are placed on radii
that satisfy condition i), then all the zeros must be placed at the
same distance $r$ away from the origin. For $r \geq R$ also condition
ii) is fulfilled. Using the symmetry we may calculate
\begin{equation*}
  \delta(B) = \frac{N r^{N-1}}{1 + r^2 + r^4 + \cdots + r^{2(N-1)}} .
\end{equation*}
\end{example}

\begin{example}
Given two arcs $\gamma_1$ and $\gamma_2$ with lengths $|\gamma_1| =
\theta$ and $|\gamma_2| = 2\pi - \theta$. We may assume that $\theta
\leq \pi$ and that the arcs lie symmetrically about the real axis. To
satisfy condition i), the zeros must be $z_1 = r_1$ and $z_2 = -r_2$
for $r_1, r_2 \in [0, 1)$, and
\begin{equation*}
  B(z) = \lambda \frac{z - r_1}{1 - r_1 z} \frac{z + r_2}{1 + r_2 z} .
\end{equation*}
As $B(\bar z) = \overline{B(z)}$ and $B(1) = \lambda$, we need to find
$r_1$ and $r_2$ such that
\begin{equation*}
  B(e^{-i \theta / 2}) = B(e^{i \theta / 2}) = -\lambda .
\end{equation*}
Solving the resulting second-order equation in $e^{i \theta / 2}$
yields
\begin{equation*}
  \cos \tfrac12 \theta = \frac{r_1 - r_2}{1 - r_1 r_2}
    \qquad \text{or} \qquad
  r_1 = \frac{r_2 + \cos \frac12 \theta}{1 + r_2 \cos \frac12 \theta}.
\end{equation*}
Because $\theta \leq \pi$ we get that $r_1 \geq r_2$, so by choosing
$r_2 \geq R$ the Blaschke product satisfies condition ii). The uniform
separation constant will be
\begin{equation*}
  \delta(B) = \frac{r_1 + r_2}{1 + r_1 r_2} .
\end{equation*}
We also remark that if $r_2 = 0$, then $r_1 = \cos \frac12
\theta$. This is the $N = 2$ version of a beautiful geometrical result
by U.\ Daepp, Gorkin and Mortini~\cite{Daepp02}.
\end{example}

\section{Constructing interpolating Blaschke products}
\label{sec:algorithm}
In this section we prove Theorem~\ref{thm:mainthm} by devising an
iterative algorithm that solves the problem. The crucial ingredients
in the algorithm are certain monotonicity relationships. To describe
these we use the harmonic measure.

Recall that for a measurable set $E \subset \partial \mathbb D$ the harmonic
measure of $E$ at a point $z \in \mathbb D$ is
\begin{equation*}
  \omega(z, E; \mathbb D) = \int_E \frac{1 - |z|^2}{|e^{i \theta} - z|^2} \,
    \frac{\mathrm d\theta}{2 \pi} .
\end{equation*}
On $\partial \mathbb D$ the derivative of the argument of a Blaschke product is
\begin{equation*}
  \tfrac{\mathrm d}{\mathrm d\theta} \bigl( \arg B(e^{i \theta}) \bigr)
    = \tfrac{\mathrm d}{\mathrm d\theta} \operatorname{Im}
      \bigl( \log B(e^{i \theta}) \bigr)
    = \sum_{n=1}^N \frac{1 - |z_n|^2}{|e^{i \theta} - z_n|^2} .
\end{equation*}
Therefore we consider the measure $\mu$ defined by
\begin{equation}
\label{eq:mu}
  \mu(E) = \sum_{n=1}^N \omega(z_n, E; \mathbb D)
    = \sum_{n=1}^N \int_E \frac{1 - |z_n|^2}{|e^{i \theta} - z_n|^2} \,
    \frac{\mathrm d\theta}{2 \pi} .
\end{equation}
With this notation our problem is to find conditions on $z_1, \ldots,
z_N$ such that $\mu(\gamma_n) = 1$ for each arc $\gamma_1, \ldots,
\gamma_N$. We first make some observations about $\mu$.

\begin{lemma}
\label{lem:mu}
The measure $\mu$ corresponding to the zeros $z_1, \ldots, z_N$ has
the following properties.
\begin{enumerate}
\item[(a)] $\mu(\partial \mathbb D) = \sum_{n=1}^N \mu(\gamma_n) = N$.
\item[(b)] $\mu(\gamma_n) \in (0, N)$.
\item[(c)] $\omega(z_n, \gamma_n; \mathbb D)$ is increasing as a function of
$|z_n|$.
\item[(d)] If $|z_n|$ is large enough, then $\omega(z_n, \gamma_m;
\mathbb D)$, $m \neq n$, is decreasing as a function of $|z_n|$.
\item[(e)] $\omega(0, \gamma_n; \mathbb D) = \frac{|\gamma_n|}{2 \pi}$,
$\lim_{|z_n| \rightarrow 1} \omega(z_n, \gamma_n; \mathbb D) = 1$ and for $n
\neq m$
\begin{equation*}
  \lim_{|z_n| \rightarrow 1} \omega(z_n, \gamma_m; \mathbb D) = 0 .
\end{equation*}
\end{enumerate}
\end{lemma}

\begin{proof}
All these properties are easy observations.  Properties (a) and (b)
follow because $\omega$ is a probability measure in the second
variable. Property (e) comes from the definition of harmonic measure,
while properties (c) and (d) follow from considerations about the
radial derivative of
\begin{equation*}
  \frac{1 - |z_n|^2}{|e^{i \theta} - z_n|^2} .
\end{equation*}
These considerations also give a sufficient condition for (d)
to hold. Namely,
\begin{equation*}
  |z_n| \geq \frac{1 - \sin \frac12 |\gamma_n|}{\cos \frac12 |\gamma_n|} .
\end{equation*}
\end{proof}

We will now describe the algorithm for constructing a sequence of
Blaschke products $(B_k)$, which converges to the Blaschke product we
seek. All the Blaschke products $B_k$ will satisfy the extra
conditions i) and ii). We denote the zeros of $B_k$ by
$z_{k,n}$. Similarily, $\mu_k$ denotes the measure defined
by~\eqref{eq:mu} corresponding to the zeros $\{z_{k,n}\}$.  To
initiate the algorithm, we calculate the mid-point of each arc and
call it $e^{i \theta_n}$, $n = 1, \ldots, N$. The zeros of the initial
Blaschke product, $B_0$, are set to be $z_{0,n} = R e^{i \theta_n}$
for $n = 1, \ldots, N$, where $R \in [0, 1)$ is chosen to be large
enough in two respects. First of all, $R$ needs to be large enough for
Lemma~\ref{lem:mu}(d) to hold. A sufficient condition for this will be
\begin{equation}
\label{eq:Rlargeenough}
  R \geq \frac{1 - \sin \frac12 L_\gamma}{\cos \frac12 L_\gamma} ,
\end{equation}
where $L_\gamma$ is the length of the shortest arc, $L_\gamma =
\min_{1 \leq n \leq N} |\gamma_n|$. Furthermore, $R$ needs to be large
enough to make $\delta(B_0) > C$. See Figure~\ref{fig:algex_step00}
for an example.

The iteration proceeds in the following manner. The Blaschke product
$B_{k+1}$ is constructed from $B_k$ by moving one of the zeros along
its corresponding radius towards the boundary $\partial \mathbb D$. To choose
which zero is moved, calculate the measures $\mu_k(\gamma_n)$ for $n =
1, \ldots, N$ corresponding to the Blaschke product $B_k$. If all
these measures are $1$, we are done. If not, identify the index of the
arc with the smallest measure and call this index $m$. The Blaschke
product $B_{k+1}$ will have the same zeros as $B_k$, except $z_m$
which is moved along its given radius so that for $B_{k+1}$ the
measure $\mu_{k+1}(\gamma_m) = 1$. See Figure~\ref{fig:algex_step01}
for an example, while Figure~\ref{fig:algorithm} gives the metacode
for the algorithm.
\begin{figure}[bt]
\begin{flushleft}
\rule{\textwidth}{0.4pt}
\textbf{Given:}\\ 
\begin{itemize}
\item $N$ disjoint arcs, $\gamma_1, \ldots, \gamma_N$.
\item The bound $C < 1$ for the separation constant.
\item The accuracy $\varepsilon > 0$.
\end{itemize}
\textbf{Algorithm:}\\ 
\begin{enumerate}
\item[1.]
Construct $B_0$.
\begin{itemize}
  \item Calculate the mid-points $e^{i \theta_n}$.
  \item Choose $R$.
  \item Set $z_{0,n} = R e^{i \theta_n}$.
  \item Set $k = 0$.
\end{itemize}
\item[2.]
Calculate the measures $\mu_k(\gamma_n)$ corresponding to the
zeros of $B_k$, and the error
\begin{equation*}
  E_k = \sum_{n=1}^N \bigl| 1 - \mu_k(\gamma_n) \bigr| .
\end{equation*}
If $E_k < \varepsilon$ then stop the algorithm.
\item[3.]
Find the index $m$ of the arc with smallest measure. (If there are
several arcs with the smallest measure, choose the index of any of
them.)
\item[4.]
Set $z_{k+1,n} = z_{k,n}$ for all $n \neq m$, and choose
$z_{k+1,m}$ such that $\mu_{k+1}(\gamma_m) = 1$.
\item[5.]
Increase $k$ by one, and return to Step 2.
\end{enumerate}
\rule{\textwidth}{0.4pt}
\end{flushleft}
\caption{Metacode for the algorithm described in Section~\ref{sec:algorithm}.}
\label{fig:algorithm}
\end{figure}

We first remark that such an iteration is always possible. If
$\mu(\gamma_n) \neq 1$ for some $n = 1, \ldots, N$, then by
Lemma~\ref{lem:mu}(a) and~\ref{lem:mu}(b) there is at least one arc
$\gamma_m$ such that $\mu(\gamma_m) < 1$. Using Lemma~\ref{lem:mu}(c)
and~\ref{lem:mu}(e) we see that there is a point on the line segment
between $|z_{k,m}| e^{i \varphi_m}$ and $e^{i \varphi_m}$ where
$z_{k+1,m}$ can be placed in order to give $\mu_{k+1}(\gamma_m) = 1$
for $B_{k+1}$.

Next, we observe that if $\mu_k(\gamma_n) \leq 1$ for some $k$, then
also $\mu_j(\gamma_n) \leq 1$ for all $j \geq k$, as the only way to
increase the measure of an arc is to move its corresponding zero, and
this will never increase the measure beyond $1$. This implies that
there is at least one arc $\gamma_l$ for which $\mu_k(\gamma_l) \geq
1$ for all $k \in \mathbb N$, and in consequence that there is at least one
zero $z_l$ which is never moved. That is, at least one zero lies at
the initial distance from the origin $R$ in all Blaschke products
$B_k$, $k \in \mathbb N$.

Finally, we comment that the sequence $(B_k)$ converges to a Blaschke
product with the desired properties. Let $B = \lim_{k \rightarrow
\infty} B_k$. Because the zeros are always moved outwards, we will
have $\delta(B) \geq \delta(B_0)$ such that a proper choice of $R$
will guarantee that $\delta(B) > C$. Thus, we only need to show that
$\mu(\gamma_n) = 1$, $n = 1, \ldots, N$, for the measure corresponding
to the zeros of $B$. For each Blaschke product $B_k$ define the error
\begin{equation*}
  E_k = \sum_{n=1}^N \bigl| 1 - \mu_k(\gamma_n) \bigr| .
\end{equation*}
Clearly, $E_k \geq 0$. Furthermore $E_{k+1} \leq E_k$ as the arcs with
$\mu_k(\gamma_n) > 1$ will contribute to a lower error, since their
measures decrease in every step. Hence $(E_k)$ is a convergent
sequence.  To see that $(E_k)$ converges to $0$ we argue that the
decrease of $E_k$ at each step is comparable to $E_k$ itself. Let
$\gamma_l$ be an arc such that $\mu_k(\gamma_l) \geq 1$ for every $k
\in \mathbb N$. Then the length of $\gamma_l$ is at least $c \frac{2\pi}N$
for some $c > 0$.  Furthermore, if $\gamma_m$ is the arc with smallest
measure at some step $k$, then $\mu_k(\gamma_m) \leq 1 - [2(N-1)]^{-1}
E_k$. When moving $z_m$, the decrease in measure from step $k$ to step
$k + 1$ is more or less evenly distributed outside the arc
$\gamma_m$. This means that the decrease $\mu_k(\gamma_l) -
\mu_{k+1}(\gamma_l)$ is at least $\tilde c [2 N (N-1)]^{-1} E_k$, and
consequently that
\begin{equation}
\label{eq:errorconvergence}
  E_k - E_{k+1} \geq \tilde c [2 N (N-1)]^{-1} E_k,
    \quad \text{for some $\tilde c > 0$}.
\end{equation}
Thus, $E_k$ converges exponentially to $0$.

\begin{example}
\label{ex:algorithm}
Assume that we are given arcs $\gamma_1, \ldots, \gamma_6$ with
lengths $\frac{\pi}{5}$, $\frac{3\pi}{5}$, $\frac{3\pi}{5}$,
$\frac{3\pi}{10}$, $\frac{\pi}{10}$ and $\frac{\pi}{5}$ respectively.
We want to find a Blaschke product $B$ with $\delta(B) > C = 0.7$,
that maps each arc onto the unit circle.
\begin{figure}[bt]
\begin{center}
\includegraphics{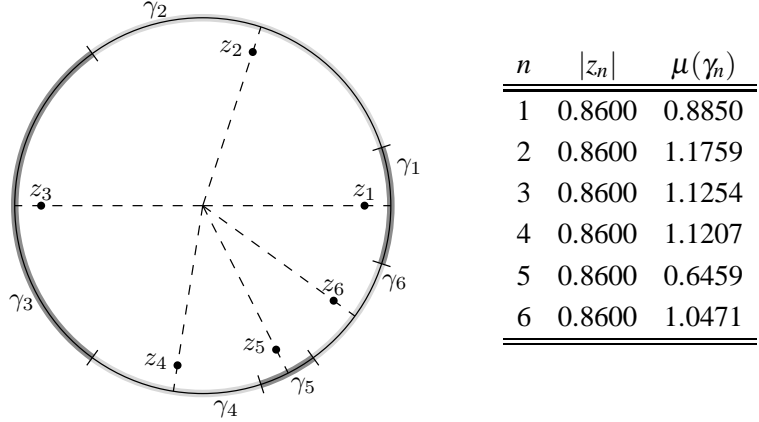}
\qquad
\renewcommand{\arraystretch}{1.2}
\raisebox{30pt}{
  \begin{tabular}[b]{cr@{.}lr@{.}l}
    $n$   & \multicolumn{2}{c}{$|z_n|$} &
      \multicolumn{2}{c}{$\mu(\gamma_n)$} \\
    \hline\hline
    1    & 0&8600 & 0&8850 \\
    2    & 0&8600 & 1&1759 \\
    3    & 0&8600 & 1&1254 \\
    4    & 0&8600 & 1&1207 \\
    5    & 0&8600 & 0&6459 \\
    6    & 0&8600 & 1&0471 \\
    \hline\hline
  \end{tabular}}
\caption{The Blaschke product $B_0$ in Example~\ref{ex:algorithm}. The error is $E_0 \approx 0.9383$.}
\label{fig:algex_step00}
\end{center}
\end{figure}

We start by constructing $B_0$, and first we calculate the mid-points.
We may assume that $\theta_1 = 0$. Then
\begin{equation*}
  \theta_2 = \tfrac{2 \pi}{5},   \quad
  \theta_3 = \pi,                \quad
  \theta_4 = \tfrac{29 \pi}{20}, \quad
  \theta_5 = \tfrac{33 \pi}{20}  \quad \text{and} \quad
  \theta_6 = \tfrac{9 \pi}{5}.
\end{equation*}
As the shortest arc is $\gamma_5$ with length $\pi / 10$, we need
\begin{equation*}
  R \geq \frac{1 - \sin (\pi / 20)}{\cos (\pi / 20)}
    \approx 0.8541
\end{equation*}
to satisfy~\eqref{eq:Rlargeenough}. Trying with $R = 0.855$, we see that
\begin{equation*}
  \delta(B_0) = \min_k \prod_{j \neq k}
    \frac{R |e^{i \theta_j} - e^{i \theta_k}|}
         {|1 - R^2 e^{\theta_j - \theta_k}|}
    \approx 0.6854 .
\end{equation*}
As this is less than $C$, we try with a bigger $R$. A new calculation
shows that $R = 0.86$ gives $\delta(B_0) \approx 0.7025 > C$, so we
use this $R$ as the initial radius.

Next, we start the iteration. First we calculate the $\mu$-measures of
the arcs for $B_0$. The result is shown in
Figure~\ref{fig:algex_step00}. We see here that $\gamma_5$ is the arc
with smallest measure. Thus, to construct $B_1$ from $B_0$ we will
move the zero $z_5$. As $\gamma_5$ has start-point $e^{i 16\pi / 10}$
and end-point $e^{i 17\pi / 10}$, we need to find conditions on
$|z_5|$ such that
\begin{equation*}
  B_1(e^{i 16\pi / 10}) = B_1(e^{i 17\pi / 10}) ,
\end{equation*}
which will imply that $\mu_1(\gamma_5) = 1$. Since we now all the
other zeros of $B_1$, this just amounts to solving a second degree
equation in $|z_5|$, and we find that $|z_5| \approx 0.9675$ is a
solution. Hence, we have found $B_1$. See
Figure~\ref{fig:algex_step01}.
\begin{figure}[bt]
\begin{center}
\includegraphics{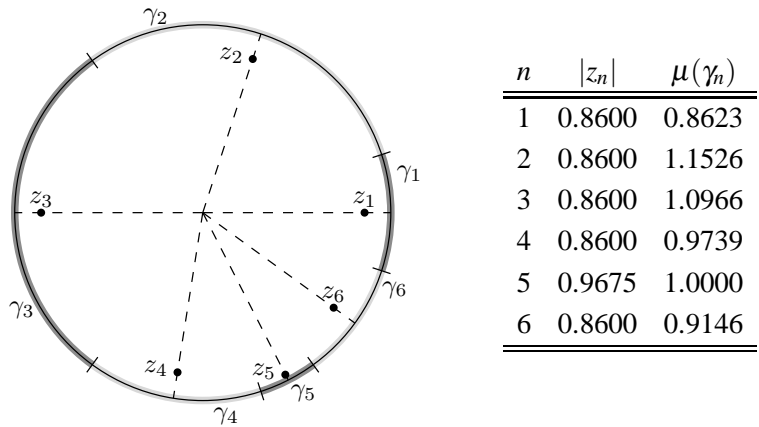}
\qquad
\renewcommand{\arraystretch}{1.2}
\raisebox{30pt}{
  \begin{tabular}[b]{cr@{.}lr@{.}l}
    $n$   & \multicolumn{2}{c}{$|z_n|$} &
      \multicolumn{2}{c}{$\mu(\gamma_n)$} \\
    \hline\hline
    1    & 0&8600 & 0&8623 \\
    2    & 0&8600 & 1&1526 \\
    3    & 0&8600 & 1&0966 \\
    4    & 0&8600 & 0&9739 \\
    5    & 0&9675 & 1&0000 \\
    6    & 0&8600 & 0&9146 \\
    \hline\hline
  \end{tabular}}
\caption{The Blaschke product $B_1$. To construct $B_1$ from $B_0$, the zero $z_5$ is moved. After 1 iteration, the error is $E_1 \approx 0.4983$.}
\label{fig:algex_step01}
\end{center}
\end{figure}

We then continue in the same manner to construct $B_2$, $B_3$, and so
on. To construct $B_2$ from $B_1$ we move the zero
$z_1$. Figure~\ref{fig:algex_step75} shows the Blaschke product
$B_{75}$. As $E_{75} \approx 1.4 \cdot 10^{-5}$ this is a quite good
approximation to the true solution $B$.
\begin{figure}[bt]
\begin{center}
\includegraphics{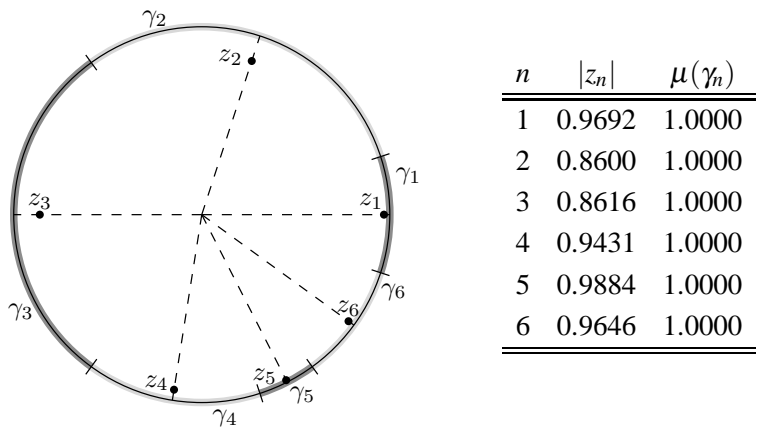}
\qquad
\renewcommand{\arraystretch}{1.2}
\raisebox{30pt}{
  \begin{tabular}[b]{cr@{.}lr@{.}l}
    $n$   & \multicolumn{2}{c}{$|z_n|$} &
      \multicolumn{2}{c}{$\mu(\gamma_n)$} \\
    \hline\hline
    1    & 0&9692 & 1&0000 \\
    2    & 0&8600 & 1&0000 \\
    3    & 0&8616 & 1&0000 \\
    4    & 0&9431 & 1&0000 \\
    5    & 0&9884 & 1&0000 \\
    6    & 0&9646 & 1&0000 \\
    \hline\hline
  \end{tabular}}
\caption{After 75 iterations, $B_k$ approximates the true solution quite well. The error $E_{75}$ is approximately $1.4 \cdot 10^{-5}$.}
\label{fig:algex_step75}
\end{center}
\end{figure}
\end{example}

\section{Variations of the result}
\label{sec:applications}

Note that there is nothing special about the mid-points of the arcs
that we chose in the extra condition i). We could let each zero move
along any radius that ends inside the corresponding arc. We do not
even need to use radii. We only need the zeros to move along curves
such that the monotonicity criteria in Lemma~\ref{lem:mu} hold. This
implies one of the strengths of our proof. Because of the flexibility
of the harmonic measure, it will apply even to more general domains
than the unit disk. The proof runs through in any domain where we can
define a measure $\mu$ that satisfies Lemma~\ref{lem:mu}. Hence, we
have the following.
\begin{corollary}
\label{cor:extension}
Let $\Gamma \subset \overline{\mathbb D}$ be a closed Jordan curve,
and let $\varphi_1, \ldots, \varphi_N$ be distinct points on
$\Gamma$. For every $C < C_\Gamma$ there is a finite Blaschke product
$B$ of degree $N$ with zeros inside $\Gamma$ such that $\delta(B) > C$
and
\begin{equation*}
  \arg B(\varphi_1) = \cdots = \arg B(\varphi_N) .
\end{equation*}
\end{corollary}
The curve $\Gamma$ puts natural restrictions on how separated the
zeros of the Blaschke product can be. This is reflected in the
constant $C_\Gamma$, which will depend on the curve $\Gamma$ and the
choice of curves that the zeros are moved along during the algorithm.

In~\cite{Gorkin05} Gorkin and Mortini proved that given a (possibly
infinite) sequence $(\lambda_k)$ of distinct points on the unit
circle, then for a sequence $(a_k) \subset \overline{\mathbb D}$ there
is an interpolating Blaschke product with radial limits $a_k$ at
$\lambda_k$ for all $k$ if and only if $(a_k)$ is bounded away from
zero. To prove this, Gorkin and Mortini needed a little more than what
is stated in Theorem~\ref{thm:mainthm}. In their paper, they showed
the following.

\begin{theorem}
\label{thm:Gorkin05}
Let $s$ be a real number satisfying $0 < s < 1$. Suppose that
$\varphi_1, \ldots, \varphi_{N+1}$ are distinct points on the unit circle
and let $\beta \in \partial \mathbb D \setminus \{ 1 \}$. Then for every $C$
with $0 < C < 1$ and $m \in \mathbb N$ there exists a Blaschke product $B$ of
degree $N$ such that
\begin{enumerate}
\item[(a)] $B(\varphi_j) = 1$ for $j = 1, \ldots, N$,
\item[(b)] $B(\varphi_{N+1}) = \beta$,
\item[(c)] $|1 - B(z)| < 2^{-m-2}$ for $|z| \leq s$,
\item[(d)] $|1 - B(r \varphi_j)| < 2^{-m}$ for $j = 1, \ldots, N$ and
$0 < r \leq 1$,
\item[(e)] $\delta(B) \geq C$.
\end{enumerate}
In addition, the following also hold:
\begin{enumerate}
\item[i)] Let $E \subset \mathbb D$. If the pseudo-hyperbolic distance
between any two distinct points in $E$ is at least $\rho$, then the
zero of $B$ closest to $\varphi_{N+1}$, can be chosen to be at
pseudo-hyperbolic distance at least $\rho/3$ to the points of $E$.
\item[ii)] It is possible to choose the zeros $z_1, \ldots, z_N$ of
$B$ so that
\begin{equation*}
  \tfrac{1 - |z_j|}{|z_j - \varphi_k|} \leq 2^{-m}
    \quad \text{for $j = 1, \ldots, N-1$ and $k = 1, \ldots, N+1$}
\end{equation*}
and such that $z_N$ is close to $\varphi_{N+1}$.
\end{enumerate}
\end{theorem}

We stated Theorem~\ref{thm:mainthm} in its simple form in order to
emphasize the ideas behind the proof. However, with some minor
modifications of our algorithm and some careful bookkeeping we can
indeed prove Theorem~\ref{thm:Gorkin05} as well.

\begin{proof}[Proof of Theorem~\ref{thm:Gorkin05}]
By identifying the points $\varphi_1, \ldots, \varphi_N$ with the
end-points of the arcs $\gamma_1, \ldots, \gamma_N$ and choosing a
proper rotation $\lambda$ we have already proved (a) and (e). To
satisfy (b) we need to add one degree of freedom to our
construction. We no longer demand that the zero $z_N$ be placed on the
radius through the mid-point of the arc $\gamma_N$. Instead we demand
that this zero be placed in the sector between $\varphi_N$ and
$\varphi_1$. This can be done using a 2-step algorithm. First we
choose some admissible radius at angle $\theta$ for the zero $z_N$ to
move along, and run the usual algorithm to construct a Blaschke
product $B_\theta$. If we find that
\begin{equation*}
  \arg B_\theta(\varphi_{N+1}) < \arg \beta
\end{equation*}
we need to choose a smaller $\theta$ and try again. Conversely, if
$\arg B_\theta(\varphi_{N+1}) > \arg \beta$ we need to choose a bigger
$\theta$. For the different angles we need to keep $R$ fixed and large
enough, where what is large enough may also depend on $\beta$. For a
properly chosen sequence of angles this process converges to a
Blaschke product which satisfies (b) in addition to (a) and (e).

Properties (c) and (d) will hold if we choose the zeros close enough
to the boundary. In this case the Blaschke product is essentially
constant and close to 1 outside the darkly shaded regions in
Figure~\ref{fig:constant}, and the disk $\{ z \colon |z| < s \}$ and
the rays $\{ z = r \varphi_j \colon j = 1, \ldots, N, 0 < r \leq 1 \}$
do not meet these regions.
\begin{figure}[bt]
\begin{center}
\includegraphics{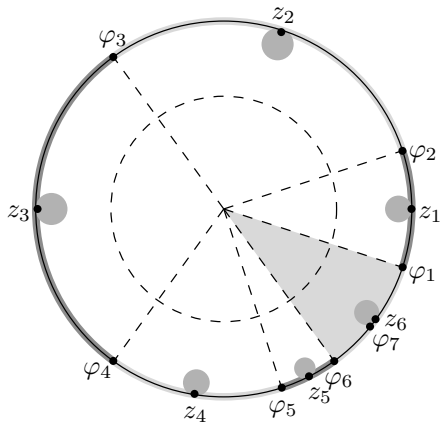}
\caption{The Blaschke product is essentially constant and close to 1 outside the darkly shaded region. The disk and the rays indicated with dashed lines do not meet these regions. The zero $z_N$ (here $N = 6$) is placed somewhere inside the lightly shaded sector.}
\label{fig:constant}
\end{center}
\end{figure}

The possible positions for the zero $z_N$ in order to fulfill (b) will
lie on a curve, parameterized by the initial radius $R$, ending in the
point $\varphi_{N+1}$. Since any point close to $\varphi_{N+1}$ lying
on the curve will yield a solution, we see that i) and the last part
of ii) holds. The first part of ii) says that the zeros $z_1, \ldots,
z_{N-1}$ can be chosen to be much closer to the boundary than to the
points $\varphi_1, \ldots, \varphi_{N+1}$, which holds true by
construction.
\end{proof}

Properties (a) and (b) in Theorem~\ref{thm:Gorkin05} give an easy way
to construct solutions to general finite interpolation problems on the
circle of the form~\eqref{eq:fip}. Namely, construct Blaschke products
$b_1, \ldots, b_N$ which each satisfy
\begin{equation*}
  b_k(\varphi_k) = \psi_k
  \qquad \text{and} \qquad
  b_k(\varphi_j) = 1 \quad \text{for $j \neq k$}
\end{equation*}
using Theorem~\ref{thm:Gorkin05}. Then $B = \prod_{k=1}^N b_k$
satisfies
\begin{equation*}
  B(\varphi_n) = \psi_n, \qquad n = 1, \ldots, N .
\end{equation*}
If we are a bit careful with the placing of the zeros of $B$, we can
also make sure that this Blaschke product has arbitrarily big
separation. However, this Blaschke product may have degree as high as
$N(N-1)$, which is far away from the optimal $N-1$ in
Theorem~\ref{thm:Jones87}. It seems plausible that it should be
possible to construct an algorithm similar to the one we have
discussed here, albeit more complicated, which can solve the general
finite interpolation problem on the circle with a Blaschke product of
degree comparable to $N$.

\bibliographystyle{amsplain}

\end{document}